\documentclass{article}
\usepackage{amsfonts,amsmath}

\newtheorem{lemma}{Lemma}
\newtheorem{theorem}{Theorem}

\newtheorem{corollary}{Corollary}
\newtheorem{proposition}{Proposition}

\def\beq{\begin{equation}}
\def\eeq{\end{equation}}

\def\R{{\mathbb R}}
\def\C{{\mathbb C}}

\def\H{{\cal H}}
\def\g{{\bf{g}}}
\def\h{{\bf{h}}}
\def\sol{\mathrm{Sol}}
\def\nil{\mathrm{Nil}}
\def\sll{\widetilde{SL(2,\R)}}
\def\D{{\cal D}}

\begin{document}

\title{Surfaces in three-dimensional Lie groups in terms of spinors}
\author{Iskander A. TAIMANOV
\thanks{Institute of Mathematics, 630090 Novosirbisk, Russia;
e-mail: taimanov@math.nsc.ru}}
\date{}
\maketitle

Recently surfaces in three-dimensional homogeneous spaces which
differ from the space forms attract a lot of attention. Mainly for
ambient spaces there are taken three-dimensional spaces with the
Thurston geometries \footnote{These are the space forms $\R^3, S^3$,
and $\H^3$; the product geometries $S^2 \times \R$ and $\H^2 \times
\R$; and three geometries modeled on the Lie groups $\nil$, $\sol$,
and $\sll$ with certain left-invariant metrics.} or simply-connected
spaces with a four-dimensional isometry group. \footnote
{\label{4isoclass} All  such spaces are locally isometric to line
bundles (with the bundle curvature $\tau$) over space forms (with
the curvature $\kappa$): for $\kappa \neq 4\tau^2$ we have the
spaces $E(\kappa,\tau)$ from the table below
$$
\begin{array}{cccc}
 & \kappa < 0 & \kappa = 0 & \kappa > 0 \\
\tau = 0 & \H^2 \times \R & \R^3 & S^2 \times \R \\
\tau \neq 0 & \sll & \nil & \mbox{Berger spheres}
\end{array}
$$
and for $\kappa = 4\tau^2$ we have spaces of constant curvature.
}

We consider the case when the ambient space is a Lie group
because it is  straightforward (see \cite{BT1}) to generalize the
Weierstrass representation of surfaces in $\R^3$ to this case. This
representation involves the Dirac operator which plays an important
role in many integrable soliton equations and has a rich and
far-developed spectral theory  \cite{T1,T2}.

In particular, we have been interested from the beginning in the following
questions:

1) it is known that certain classes of surfaces in the space forms
are  described by some integrable systems
(for instance, constant mean curvature tori).

{\it How such surfaces are described in new geometries?}

{\it If these surfaces are described by some integrable systems
how these systems obtained from the old
ones and how the curvature of the ambient space contribute to the deformation
of an integrable system?}

2) it is known that some spectral data of the Dirac operator coming
in the Weierstrass representation of a surface in $\R^3$ have
geometrical meanings and, in particular, the Willmore functional
serves as an example \cite{T1,T2}.

{\it What mean these spectral data for surfaces in other
ambient spaces (Lie groups)?}

We discuss some partial answers to these questions in \S\S 4 and 5.

We also would like to mention that the choice of Lie groups as the
ambient spaces is not very restrictive since it covers all spaces
$E(\kappa,\tau)$ but $S^2 \times \R$ and the Thurston geometries
again except $S^2 \times \R$ (see remarks on page \pageref{lieiso}).

{\bf 1. The Weierstrass representation of surfaces in $\R^3$
and the Willmore functional}

The original Weierstrass representation of minimal surfaces in
$\R^3$ may be considered as an integrable system in geometry because
it gives an explicit formula for a general solution to the minimal
surface equation in $\R^3$ in terms of a pair of arbitrary
holomorphic functions. It is as follows. Let $z \in D \subset \C$
and, for simplicity, assume that a domain $D$ is simply connected.
Let $f$ and $g$ be holomorphic functions on $D$. Then the
Weierstrass (--Enneper) formulas
$$
x^1(z,\bar{z})
= x^1_0 +
\frac{i}{2} \int \left[(f^2+g^2)dz - (\bar{f}^2+\bar{g}^2)d\bar{z})\right],
$$
\beq
\label{weierstrass}
x^2 = x^2_0 + \frac{1}{2} \int
\left[(g^2-f^2)dz + (\bar{g}^2 - \bar{f}^2)d\bar{z})\right],
\eeq
$$
x^3 = x^3_0 + \int (fg dz + \bar{f}\bar{g} d\bar{z})
$$
define a minimal surface in $\R^3$. Here the integrals defining
$x(P)$, the image of $P \in D$, are taken along a path $\gamma
\subset D$ from the point $P_0$ such that $x(P_0)=x_0$ to $P$. Since
the integrands are closed forms this is independent on the choice of
$\gamma$. The induced metric takes the form $(|f|^2+|g|^2)^2 dz
d\bar{z}$ and therefore $z$ is a conformal parameter on the surface

In fact, the condition that $z$ is a conformal parameter is written as
$$
\left(\frac{\partial x^1}{\partial z}\right)^2 +
\left(\frac{\partial x^2}{\partial z}\right)^2 +
\left(\frac{\partial x^3}{\partial z}\right)^2 = 0,
$$
i.e., $(r_u,r_u)=(r_v,r_v), (r_u,r_v)=0$ where $u$ and $v$ are
the isothermic coordinates such that
$z = u+iv$, and $r_u = 2 \mathrm{Re} \frac{\partial x}{\partial z}$ and
$r_v = -2 \mathrm{Im} \frac{\partial x}{\partial z}$ are
the corresponding tangent vectors to the surface.
The quadric
$$
Q = \{y_1^2 + y_2^2 + y_3^2 = 0\} \subset  \C P^2
$$
gives a one-to-one parametrization of oriented two-planes in $\R^3$
by corresponding to every plane its homogeneous coordinates $((\xi^1
- i \eta^1): (\xi^2 - i \eta^2) : (\xi^3 - i \eta^3))$ where $(\xi,
\eta)$ is a positively oriented basis for the plane such that $|\xi|
= |\eta|$ and $\xi$ is orthogonal to $\eta$. Due to the homogeneity
of coordinates in $\C P^2$ this mapping is correctly defined, i.e.,
is independent on the choice of a basis $(\xi,\eta)$. Hence the
mapping
$$
P \to \left(\frac{\partial x^1(P)}{\partial z} :
\frac{\partial x^2(P)}{\partial z} :
\frac{\partial x^3(P)}{\partial z}\right) \in Q
$$
is the Gauss map of the surface. The quadric $Q$, the Grassmannian
of oriented two-planes in $\R^3$, admits a natural rational
parametrization: \beq \label{factorization} (f:g) \to
\left(\frac{i}{2}(f^2+g^2): \frac{1}{2}(g^2-f^2):fg\right). \eeq

From this interpretation of the Gauss map it is clear that

\begin{itemize}
\item
any surface, not only minimal,
is defined by the Weierstrass formulas for the factorization $(f,g)$ of the
Gauss map.
\end{itemize}

The Gauss--Codazzi equations written in terms of $(f,g)$
distinguish mappings $D \stackrel{(f:g)}{\longrightarrow} Q$
which are the Gauss maps of surfaces. It is straightforward to compute that
these equations take the form
$$
\D \psi = 0
$$
where $\D$ is the Dirac operator
\beq
\label{dirac}
\D = \left(\begin{array}{cc} 0 & \partial \\ -\bar{\partial} & 0
\end{array} \right) +
\left(\begin{array}{cc} U & 0 \\ 0 & V
\end{array} \right)
\eeq
and
$$
\psi = \left(\begin{array}{c} \psi_1 \\ \psi_2 \end{array}\right) \left(\begin{array}{c} f \\ \bar{g} \end{array}\right).
$$
For surfaces in $\R^3$ the potentials $U$ and $V$ and the induced metric
are
\beq
\label{data}
U = V = \frac{He^\alpha}{2}, \ \ \
e^{2\alpha} dzd\bar{z} = (|\psi_1|^2+|\psi_2|^2)^2dzd\bar{z}.
\eeq
We conclude that
\begin{itemize}
\item
a general surface in $\R^3$ is represented by the Weierstrass
formulas (\ref{weierstrass}) for some solution to the Dirac equation
with the potentials (\ref{data}) and the inverse is also true: any
solution to the Dirac equation with real-valued potentials $U=V$
defines via (\ref{weierstrass}) a surface in $\R^3$ with the mean
curvature and the induced metric given by (\ref{data}).
\end{itemize}

This representation has some prehistory for which we refer to
\cite{T2} however for $U \neq 0$ the formulas in terms of the Dirac
operator first appeared in \cite{Kon} where they were introduced for
inducing surfaces admitting certain soliton deformations. This
operator has a rich spectral theory and, in particular, we started
in \cite{T1} to study  possible relations between the spectral
properties of $\D$ and the geometry of the corresponding surfaces.
In particular, it appears that for a closed oriented surface $M
\subset \R^3$ the integral \beq \label{energy} E(M) = \int_M UV dx
dy \eeq is one-fourth of the Willmore functional \beq
\label{willmore} {\cal W}(M) = \int_M H^2 d\mu \eeq where $d\mu$ is
the induced measure on $M$. The Willmore functional is the basic
functional in the conformal surface geometry, and the integral
(\ref{energy}) is an important spectral quantity of the Dirac
operator $\D$.

The Willmore conjecture states that ${\cal W}$ attains its minima
for tori which is equal to $2\pi^2$ on the Clifford torus and its
images under conformal transformations of $\bar{\R^3}$. The
existence of the lower bounds for ${\cal W}$ on closed surfaces is
explained by the Weierstrass representation as follows:

\begin{itemize}
\item
there are no compact minimal surfaces without boundary in $\R^3$.
We have to perturb the potential $U$ from the zero level to achieve
compact surfaces and the threshold for the $L_2$-norm of $U$ at which compact
surfaces appear gives this minimum level. For surfaces in $\R^3$ we have
$U=\bar{U}=V$, the energy (\ref{energy}) is the squared $L_2$-norm of $U$ and
it is also one-fourth of ${\cal W}$.
\end{itemize}

We propose an approach to the Willmore conjecture based on the
spectral properties of the corresponding double-periodic (for tori)
Dirac operator. Several attempts to realize this approach led to
interesting results however the conjecture stays open until
recently. We refer for the survey of the Willmore conjecture and the
spectral approach to its study to \cite{T2}.

The classical Weierstrass representation for minimal surfaces
corresponds to the case $U = 0$ and it enables us to consider the
minimal surface equation in $\R^3$ as an integrable system. The
integrability property resolves the local theory and does not help
straightforwardly in answering questions on the global behavior of
surfaces. The global theory needs an additional technique concerning
global solutions to the integrable system (in the case of minimal
surfaces, holomorphic functions).

{\bf 2. The Weierstrass representation of surfaces in
three-dimensional Lie groups \cite{BT1}}

To generalize the Weierstrass representation for the case when
the ambient space is a three-dimensional Lie group $G$ with
a left-invariant metric \cite{BT1}
we have to replace  $\frac{\partial x}{\partial z} \in \C^3$
by the element of the complexified Lie algebra:
$$
\frac{\partial}{\partial z} \in \C^3 \longrightarrow
\Psi = f^{-1} \frac{\partial f}{\partial z} \in \g \otimes \C
$$
where
$$
f: M \to G
$$
is an immersion of a surface and $z$ is a conformal parameter on
$M$. In terms of $\Psi$ and $\Psi^\ast = f^{-1} f_{\bar{z}} \bar{\Psi}$ the derivational equations take the form
$$
\partial\Psi^\ast - \bar{\partial}\Psi + \nabla_{\Psi}\Psi^\ast -
\nabla_{\Psi^\ast}\Psi = 0,
$$
$$
\partial\Psi^\ast + \bar{\partial}\Psi + \nabla_{\Psi}\Psi^\ast +
\nabla_{\Psi^\ast}\Psi = e^{2\alpha} H f^{-1}(N)
$$
where the Levi-Civita connection on $G$ is linearly expanded onto
complex-valued vectors $\Psi$ and $\Psi^\ast$, $N$ is the unit
normal vector field to $M$ and $e^{2\alpha} dzd\bar{z}$ is the
induced metric. Originally these equations were first derived for
minimal surfaces in \cite{Hitchin}.

Given an orthonormal basis $e_1,e_2,e_3$ for $\g$, we expand
$\Psi$ in this basis
$$
\Psi = Z_1 e_1 + Z_2 e_2 + Z_3 e_3.
$$
The conformality condition again takes the form
$$
Z_1^2 + Z_2^2 + Z_3^2 = 0.
$$
Let us use the same factorization of $Z: M \to Q$ as in the
Euclidean case:
$$
Z_1 = \frac{i}{2} ( \bar{\psi}_2^2 + \psi_1^2), \ \ \
Z_2 = \frac{1}{2} ( \bar{\psi}_2^2 - \psi_1^2), \ \ \
Z_3 = \psi_1 \bar{\psi}_2.
$$
The derivational equations take the form of the Dirac equation
$$
\D \psi = 0
$$
and the induced metric is again equal to
$$
e^{2\alpha} d zd\bar{z} = (|\psi_1|^2+|\psi_2|^2)^2 d zd\bar{z}.
$$
Therewith we call $\psi$ a generating spinor of a surface.

In difference with the Euclidean case,
the potentials $U$ and
$V$ are not always real-valued and do not always coincide.

\begin{theorem}
[\cite{BT1}] \label{theorem1} The potentials of the Weierstrass
representation of surfaces in the Lie groups $SU(2), \nil, \sll$,
and $\sol$, endowed with the Thurston geometries, are as follows:

\begin{enumerate}
\item
$G = SU(2)$:
$$
U = \bar{V} = \frac{1}{2}(H-i)e^\alpha;
$$

\item
$G = \nil$:
$$
U=V = \frac{He^\alpha}{2} +
\frac{i}{4} (|\psi_2|^2 - |\psi_1|^2);
$$

\item
$G= \sll$:
$$
U = \frac{He^\alpha}{2} + i\left(\frac{1}{2}|\psi_1|^2 -
\frac{3}{4}|\psi_2|^2\right),\ \
V = \frac{He^\alpha}{2} + i\left(\frac{3}{4}|\psi_1|^2 -
\frac{1}{2}|\psi_2|^2\right);
$$

\item
$G=\sol$:
\footnote{Here we correct the sign of the second term in the expression
for $U$ miscalculated in \cite{BT1}.}
$$
U = \frac{He^\alpha}{2} -
\frac{1}{2}\bar{\psi}_2^2 \frac{\bar{\psi}_1}{\psi_1}, \ \ \
V = \frac{1}{2}He^\alpha +
\frac{1}{2}\bar{\psi}_1^2 \frac{\bar{\psi}_2}{\psi_2}.
$$
\end{enumerate}
\end{theorem}

These potentials are written with respect to certain choices of
orthogonal bases for
$\g$ which are as follows:

a) $\sol$ admits a natural splitting
$$
1 \to \R^2 \to \sol \to \R
$$
which induces the submersion $\sol \to \R = \sol/\R^2$ whose leaves
are minimal surfaces. We put $e_3$ to be the pullback of the unit
vector on $\R$. Hence, $Z_3 = \psi_1 \bar{\psi}_2 =0$ if the tangent
plane to a surface is tangent to a minimal leave. For a surface in
$\sol$ the Dirac equation is correctly defined only in domain
$D=\{Z_3 \neq 0\}$. It is natural to assume that $U=V=0$ outside
$D$. Then the Dirac equations hold everywhere outside $\partial D$,
the boundary of $D$, at which $\frac{\bar{\psi}_1}{\psi_1}$ and
$\frac{\bar{\psi}_2}{\psi_2}$ may have indeterminancies;

b) for $\nil$ and $\sll$ we assume that $e_3$ is directed along the axis of
isometry rotation.
Both these groups admit four-dimensional isometry groups and such an axis is
uniquely defined everywhere.

These Dirac equations differ from their Euclidean analog in several
aspects:

a) there are constraints which relate solutions $\psi$
corresponding to surfaces with potentials. In the Euclidean case
any solution corresponds to a surface. This demonstrates the
absence of dilations in these Lie groups;

b) the reconstruction of the surface $f: M \to G$ from $\psi$ needs solving
the linear equation
$$
f_z = f \Psi.
$$
In the Euclidean case a solution to this equation is given
by (\ref{weierstrass});

c) solutions to these Dirac equation does not admit the quaternion symmetry,
i.e., if
$\left(\begin{array}{c} \psi_1 \\ \psi_2 \end{array}\right)$ satisfies
$\D\psi=0$
then in general $\psi^\ast = \left(\begin{array}{c} -\bar{\psi}_2 \\
\bar{\psi}_1 \end{array}\right)$ does not meet this equation. This hinders to
use the Dirac equation for interpretting
surfaces as holomorphic sections of certain line bundles
and applying some ideas of algebraic geometry as it is done for surfaces in
$\R^3$ and $\R^4$ in \cite{FLPP}.

\begin{corollary}
\footnote{We skip here the well-studied case of minimal surfaces
in the unit three-sphere $SU(2)$.}
The generating spinors of minimal surfaces
in the Lie groups $\nil, \sll$, and $\sol$
are given by the following equations:

\begin{enumerate}
\item
$G= \nil$:
$$
\bar{\partial} \psi_1 = \frac{i}{4}(|\psi_2|^2 - |\psi_1|^2)\psi_1, \ \ \
\partial \psi_2 = -\frac{i}{4}(|\psi_2|^2 - |\psi_1|^2)\psi_2;
$$

\item
$G=\sll$:
$$
\bar{\partial} \psi_1 i\left(\frac{3}{4}|\psi_1|^2-\frac{1}{2}|\psi_2|^2\right)
\psi_2,
\ \ \
\partial \psi_2 -i\left(\frac{1}{2}|\psi_1|^2-\frac{3}{4}|\psi_2|^2\right)\psi_1;
$$

\item
$G=\sol$:
$$
\bar{\partial}\psi_1 = \frac{1}{2} \bar{\psi}_1^2 \bar{\psi}_2, \ \ \
\partial \psi_2 = \frac{1}{2} \bar{\psi}_1 \bar{\psi}_2^2.
$$
\end{enumerate}
\end{corollary}

In other terms the Weierstrass type representations for minimal surfaces in
$\nil$ and $\sol$ were derived in \cite{I1,I2}.

We remark that Friedrich showed that the $\psi$-spinor for surfaces
in $\R^3$ may be interpreted as the restriction of the parallel
spinor field in $\R^3$ onto the surface \cite{Friedrich}. Later a
similar description of such representations for surfaces in $S^3$
and $\H^3$ was derived in \cite{Morel} and very recently the same
was done for surfaces in the spaces with a four-dimensional isometry
group \cite{Roth} (this paper uses description of immersions in
other terms obtained in \cite{Daniel}). In the first case the
parallel spinor field is replaced by real and imaginary Killing
fields and and in the second case it is replaced by certain
generalized Killing spinor fields.

{\bf 3. Surfaces in general Lie groups and families of Lie groups.}

The Weierstrass representation method admits us to write such
representations straightforwardly for a general Lie group and even
to consider surfaces in families of Lie groups. We demonstrate that
for a certain family which includes some well-known spaces.

Let us remind {\sc the Bianci classification} of real three-dimensional
Lie algebras.

For such an algebra $\g$ there is a basis $e_1,e_2,e_3$ such that
the commutation relations takes the form
$$
[e_1,e_2]
= ae_2 + b^{(3)}e_3, \ \ [e_1,e_3] = ae_3 - b^{(2)}e_2, \ \
[e_2,e_3] = b^{(1)}e_1
$$
with $ab^{(1)} = 0$, hence the Lie
algebra is included in the following table
$$
\begin{array}{cccccccccc}
\mathrm{Type} & a & b^{(1)} & b^{(2)} & b^{(3)} & \mathrm{Type} & a & b^{(1)}
& b^{(2)} & b^{(3)} \\
\mathrm{I} & 0 & 0 & 0 & 0 & \mathrm{VI}_0 & 0 & 1 & -1 & 0 \\
\mathrm{II} & 0 & 1 & 0 & 0 & \mathrm{VI}_a, 0<a<\infty, a \neq 1 & a & 0 & 1
& -1 \\
\mathrm{III} & 1 & 0 & 1 & -1 &\mathrm{VII}_0 & 0 & 1 & 1 & 0 \\
\mathrm{IV} & 1 & 0 & 0 & 1 & \mathrm{VII}_a, a>0 & a & 0 & 1 & 1\\
\mathrm{V} & 1 & 0 & 0 & 0 & \mathrm{VIII} & 0 & 1 & 1 & -1 \\
& & & & & \mathrm{IX} & 0 & 1 & 1 & 1
\end{array}
$$
and algebras corresponding to different entries of this table are
pairwise nonisomorphic.

The simply-connected Lie groups with Lie
algebras of types $\mathrm{I}$--$\mathrm{VII}$ have the form
\beq
\label{ext}
1 \to \R^2 = H \to G \to G/H = \R \to 1
\eeq
and such an extension is uniquely defined by the action
$$
\mathrm{Ad}_z X = zXz^{-1} = e^{Az} X, \ \ z \in G/H,
\ X =  \left(\begin{array}{c} x \\ y \end{array}\right) \in H,
\ A \in gl(2,\R).
$$
In terms of Lie algebras we have
$$
\mathrm{ad}_\eta \xi = [\eta,\xi] = A\xi
$$
where $\eta$ and the Lie algebra $\h$ of $H$ span $\g$ and $\xi
\in \h$. The matrices $A$ and $\lambda BAB^{-1}$, $\lambda =
\mathrm{const}\neq 0$, define isomorphic extensions.

We have

$\mathrm{I}$: $G = \R^3, A=0$.

$\mathrm{II}$: $G = \nil$, the nilpotent group, $A = \left(\begin{array}{cc} 0 & 1\\
0 & 0 \end{array}\right)$.

$\mathrm{III}$: $G = \R \times A(1)$, where $A(1)$ is the group of
all affine transformations of $\R^1$;
$A =\left(\begin{array}{cc} 1 & 0\\
0 & 0 \end{array}\right)$.

$\mathrm{IV}$: $A = \left(\begin{array}{cc} 1 & 1\\
0 & 1 \end{array}\right)$.

$\mathrm{V}$: $G=A(2)$, the group formed by three-dimensional
affine transformations of the form \beq \label{hyperbolic}
\left(\begin{array}{cc} e^t \cdot I_2 & s \\ 0 & 1
\end{array}\right), \ \ \ t \in \R, s \in \R^2; \eeq
$A = \left(\begin{array}{cc} 1 & 0\\
0 & 1 \end{array}\right)$, $I_2$ is the unit $(2 \times
2)$-matrix.

$\mathrm{VI}_0$: $G = \sol$, the solvable group;
$A = \left(\begin{array}{cc} 1 & 0 \\ 0 & -1
\end{array}\right)$.

$\mathrm{VI}_a, a \neq 0$:
$A = \left(\begin{array}{cc} a & -1 \\ -1 & a \end{array}\right)$,
the eigenvalues $\lambda_{1,2}$ of $A$ are $\lambda_{1,2} = a \pm 1$.

$\mathrm{VII}_0$: $G = E(2)$, the group of all
isometries of $\R^2$; $A = \left(\begin{array}{cc} 0 & -1 \\
1 & 0
\end{array}\right)$.

$\mathrm{VII}_a, a \neq 0$:
$A = \left(\begin{array}{cc} a & 1 \\ -1 & a \end{array}\right)$,
the eigenvalues of $A$ are $\lambda_{1,2} = a \pm i$.

The algebras
of types $\mathrm{VIII}$ and $\mathrm{IX}$ do not contain two-dimensional
commutative subalgebras and hence does not admit the representation
(\ref{ext}). We have

$\mathrm{VIII}$: $G = \widetilde{SL(2,\R)}$, the universal cover of $SL(2,\R)$,
which is also locally isomorphic to $SO(1,2)$ and $SU(1,1)$.

$\mathrm{IX}$: $G=SU(2) = \widetilde{SO(3)}$.

A left-invariant metric on a Lie group $G$ is uniquely
defined by its value at the unit of $G$, i.e. by an inner product on
the Lie algebra $\g$. Given an orthonormal basis
$e_1,\dots,e_n$ for $\g$:
$\langle e_i,e_j \rangle = \delta_{ij}, \ \ \ i,j =1, \dots,n$,
we denote by the same symbols the corresponding left-invariant vector fields.
The Levi-Civita connection is given by the following formulas:
$$
\nabla_{e_k} e_j = \Gamma^i_{jk} e_i, \ \
\Gamma^i_{jk} = \frac{1}{2} \left(c^i_{kj} + c^j_{ik} +
c^k_{ij}\right), \ \ \
[e_i, e_j] = c^k_{ij} e_k.
$$

Let us denote by $H_n$ the group of all $n$-dimensional affine
transformations of the form (\ref{hyperbolic}) with $s \in
\R^{n-1}$. By simple computations we obtain

\begin{proposition}
\footnote{Recently we have known that such a representation of the
hyperbolic three-space was used by Kokubu for deriving the
Weierstrass representation of minimal surfaces in $\H^3$
\cite{Kokubu}.} Let us endow the group $H_n$ by the left-invariant
metric for which $e_1 = \frac{\partial}{\partial t}, e_2 =
\frac{\partial}{\partial s^1}, \dots, e_n =
\frac{\partial}{\partial s^{n-1}}$ for the orthonormal basis in
$\g$. Then $H_n$ is isometric to the $n$-dimensional hyperbolic
space $\H^n$.
\end{proposition}

\begin{corollary}
The group of type $\mathrm{III}$ with a certain left-invariant metric is
isometric to $\H^2 \times \R$.
\end{corollary}

\begin{corollary}
There is a left-invariant metric on the group of type $\mathrm{V}$
such that such a Riemannian manifold is isometric to  $\H^n$.

$H_n$ acts isometrically by left translations on
$\H^n = \{(x,y), x \in \R^{n-1}, y \in \R, y >0\}$ with the metric
$ \frac{dx^2 + dy^2}{y^2}$ as follows: $(x,y) \to (e^tx+s,e^ty).$
\end{corollary}

We see that

\begin{itemize}
\item
\label{lieiso}
all simply-connected homogeneous three-spaces with a four-dimensional iso\-metry
group except $S^2 \times \R$
are isometric to Lie groups with left-invariant metrics

\item
all Thurston geometries except $S^2 \times \R$ are modeled by Lie groups
with left-invariant metrics.
\end{itemize}

Let us consider the $\mu$-parameter family $G_\mu$ of Lie groups of type
(\ref{ext}) for which
$$
A_\mu = \left(\begin{array}{cc}
\mu & 0 \\ 0 & 1 \end{array}\right).
$$
For $-1 \leq \mu \leq 1$ these groups are pairwise nonisomorphic and
as follows:

$\mu = -1$: \ \ $\sol$, i.e. of the type $\mathrm{VI}_0$;

$-1 < \mu < 0$: \ \ $\mathrm{VI}_a, \ 0 < a <1, \ \mu = \frac{a-1}{a+1}$;

$\mu = 0$: \ \ $\mathrm{III}$;

$0 < \mu < 1$: \ \ $\mathrm{VI}_a, \ 1 < a < \infty, \ \mu = \frac{a-1}{a+1}$;

$\mu =1$: \ \ $\mathrm{V}$.

Let us take the orthonormal basis $e_1,e_2,e_3$ such that
$$
[e_1,e_2] = 0, \ \ [e_3,e_1] = \mu e_1, \ \ [e_3,e_2] = e_2.
$$
For the corresponding left-invariant metrics we have
$$
G_{-1} = \sol, \ \ G_0 = \H^2 \times \R,\ \ G_1 = \H^3.
$$

\begin{proposition}
The potentials of the Weierstrass representation for surfaces in $G_\mu$
are as follows:
$$
U_\mu = \frac{H}{2} e^{\alpha} + \frac{\mu +1}{4} |\psi_1|^2 +
\frac{\mu -1}{4} \frac{\bar{\psi}_2 ^2
\bar{\psi}_1}{\psi_1},
$$
$$
 V_\mu = \frac{H}{2} e^{\alpha} -
\frac{\mu +1}{4} |\psi_2|^2 - \frac{\mu -1}{4}
\frac{\bar{\psi}_1 ^2 \bar{\psi}_2}{\psi_2}.
$$

The generating spinor $\psi$ of a minimal surface in $G_\mu$ meets
the equations
$$
\bar{\partial}\psi_1 = -\frac{\mu + 1}{4}\psi_2 ^2 \bar{\psi}_2 -
\frac{\mu -1}{4}\bar{\psi}_1 ^2 \bar{\psi} _2,
$$
$$
\partial \psi_2 = -\frac{\mu + 1}{4}\psi_1 ^2 \bar{\psi}_1 -
\frac{\mu-1}{4}\bar{\psi}_2 ^2 \bar{\psi} _1.
$$
\end{proposition}

In early 1900s for proving the existence of three closed
nonselfintersecting geodesics on a two-sphere with a general
metric, Poincare proposed to take an analytical $\mu$-parameter
family of metrics which joins the metric on the ellipsoid with
three different axes and the given metric and then to consider the
analytical continuation in $\mu$ of the plane sections of the
ellipsoid. This program was not realized however it led to some
interesting results on perturbations of closed geodesics under
deformations of metrics.

It also would be interesting to study the $\mu$-deformations of
integrable surfaces in $G_\mu$. Probably that could help to

{\sl extend some global results on well-studied minimal or, more general, constant
mean curvature surfaces
in $G_1 = \H^3$ to such surfaces in $\sol$.}

{\bf 4. Constant mean curvature (CMC) surfaces in Lie groups}

The second fundamental form of a surface in $\R^3$ is uniquely determined by
the mean curvature $H$ and the Hopf quadratic differential
$$
A dz^2 = (x_{zz},N)dz^2,
$$
where $x_{zz} = \frac{\partial^2 x}{\partial z^2}$ and $N$ is
the unit normal vector field.
We have
$$
|A|^2 = \frac{(\kappa_1 - \kappa_2)^2e^{4\alpha}}{16}
$$
where $\kappa_1$ and $\kappa_2$ are the principal curvatures.
In terms of $\psi$ this differential takes the form
$$
A = \bar{\psi}_2 \partial \psi_1 - \psi_2 \partial \bar{\psi}_2.
$$
The Gauss--Codazzi equations are
\beq
\label{gauss}
\alpha_{z\bar{z}} + U^2 - |A|^2 e^{-2\alpha} = 0
\eeq
which is the
Gauss formula for the curvature in terms of the metric and
$$
A_{\bar{z}} = (U_z - \alpha_z U)e^\alpha
$$
which implies that $A$
is holomorphic if and only if $H = \mathrm{const}$.
\footnote{The analogous results were established by
Hopf also for surfaces in other
space forms, $S^3$ and $\H^3$.}

Since the only holomorphic quadratic differential on a sphere
vanishes everywhere, any CMC sphere in $\R^3$ is umbilic, i.e.,
$\kappa_1=\kappa_2$ everywhere, and it is easily to derive that any
closed umbilic surface is a round sphere. For tori the holomorphic
quadratic differentials are constant and, since there are no umbilic
tori, the Hopf differential of a CMC torus equals $\mathrm{const}\cdot
dz^2 \neq 0$. By a dilation any CMC torus is transformed into the
torus with $H=1$ and then by rescaling a conformal parameter we may
achieve $A = \frac{1}{2}$. Then (\ref{gauss}) takes the form
\beq
\label{sinhgordon}
u_{z\bar{z}} + \sinh u =0, \ \ \ u=2\alpha,
\eeq
which is the integrable elliptic $\sinh$-Gordon equation (see the
classification of such tori based on this integrable system in
\cite{PinkallSterling}).

Recently such an approach was extended for studying CMC surfaces
in other ambient spaces. The breakthrough point was a result of
Abresch and Rosenberg who proved that

\begin{itemize}
\item
{\sl there is a generalized Hopf differential $A_{\mathrm{AR}}dz^2$
which is defined on any surface in $S^2 \times \R$ or $\H^2 \times
\R$ such that for CMC surfaces $A_{\mathrm{AR}}$ is holomorphic}
\end{itemize}

\noindent by deriving the explicit formula for this differential
\cite{AR1}. This differential vanishes identically on a CMC sphere
and they are shown that if the equations $H = \mathrm{const}$ and
$A_{\mathrm{AR}} = 0$ are satisfied on a closed surface $M$ then
$M$ is a sphere of revolution which implies that

\begin{itemize}
\item
{\sl every CMC sphere in $S^2 \times \R$ or $\H^2 \times \R$ is a sphere
of revolution.}
\end{itemize}

Later they extended that for surfaces in other homogeneous
manifolds with a four-dimensional isometry group \cite{AR2}.
Moreover Abresch announced that

\begin{itemize}
\item
{\sl only the spaces $E(\kappa,\tau)$ admit generalized Hopf differentials
which are holomorphic on CMC surfaces.}
\end{itemize}

The mashinery of the Weierstrass representation admits us to derive very
easily such differentials for surfaces in $\nil$ and $\sll$ and moreover
to study (the first time) the following problem:

{\sl When the holomorphicity of the generalized Hopf differential
implies that the surface has constant mean curvature?}

It appeared that although for $\nil$ the answer is positive as for
space forms in general, there are non-CMC surfaces with
holomorphic generalized Hopf differential (see \cite{FM} and
below).

We have

\begin{theorem}
[\cite{BT1}]
Let us denote by $Adz^2 = (\nabla_{f_z}f_z,N)dz^2$ the Hopf differential
of a surface $f: M \to G$. Then

\begin{enumerate}
\item
for $G= \nil$ the quadratic differential
\beq
\label{tildea}
\widetilde{A} dz^2 = \left (A + \frac{{Z_3}^2}{2H+i} \right)dz^2
\eeq
is holomorphic on a surface if and only if the surface has constant
mean curvature;

\item
for $G= \sll$ the quadratic differential
$$
\widetilde{A} dz^2 = \left(A + \frac{5}{2(H-i)}Z^2_3\right) dz^2
$$
is holomorphic on constant mean curvature surfaces.
\end{enumerate}
\end{theorem}

The original Abresch--Rosenberg differential $A_{\mathrm{AR}}$
derived in \cite{AR1,AR2} is slightly
different from ours:
$$
A_{\mathrm{AR}} = (H+i\tau) \widetilde{A}
$$
where $\tau$ is the bundle curvature (see footnote on page
\pageref{4isoclass}).
These differentials behave differently for non-CMC surfaces.
Fernandez and Mira \cite{FM} showed how the definition of
$\widetilde{A}$ is extended for other spaces $E(\kappa,\tau)$
and proved that

\begin{itemize}
\item
{\sl a compact surface $M \subset E(\kappa,\tau)$ with holomorphic differential $\widetilde{A}$
(if $\tau \neq 0$ we assume that $M$ is not a torus)
is a CMC surface;}

\item
{\sl in $\H^2 \times \R$ and $\sll$ all surfaces with holomorphic
differential $\widetilde{A}$ are CMC-surfaces or some non-compact
surfaces whose complete description is given in \cite{FM};}

\item
{\sl there are non-compact rotationally-invariant non-CMC surfaces
with holomorphic differential $A_{\mathrm{AR}}$ in $S^2 \times \R$
and $\H^2 \times \R$ however it is still unclear are there non-CMC surfaces
with holomorphic differential $\widetilde{A}$ in such ambient
spaces.}
\end{itemize}

As we see above CMC-tori in $\R^3$ are described by the
elliptic $\sinh$-Gordon equation. By \cite{BT1,FM}, in
the spaces $E(\kappa,\tau)$ except probably some Berger spheres CMC
tori are exactly the tori with holomorphic differential
$\widetilde{A}$. It appeared that for surfaces in $\nil$ the
holomorphicity of $\widetilde{A}$ again leads to the elliptic $\sinh$-Gordon
equation but for other quantities.

\begin{theorem}
[Berdinsky]
For a certain choice of a conformal parameter the potential $U=V$ of
the Weierstrass representation of a CMC torus has to meet the equation
\beq
\label{berdinsky}
v_{z\bar{z}} + 2\sinh 2v = 0
\eeq
where $v = \log U$.
\end{theorem}

First we prove the following

\begin{lemma}
[Berdinsky]
In terms of $\psi$ and of the differential
$$
B =  \frac{1}{4} (2H+i) \widetilde{A}
$$
the derivational equations for surfaces in $\nil$ are written as follows
\beq
\label{dernil1}
\partial \left(
\begin{array}{c}
\psi_1 \\
\psi_2
\end{array}
\right)
 \left(
\begin{array}{cc}
v_z - \frac{1}{2} H_z e^{-v} e^{\alpha} & Be^{-v}\\
-e^v & 0
\end{array}
\right)
\left(
\begin{array}{c}
\psi_1 \\
\psi_2
\end{array}
\right),
\eeq
\beq
\label{dernil2}
\bar{\partial} \left(
\begin{array}{c}
\psi_1 \\
\psi_2
\end{array}
\right)
 \left(
\begin{array}{cc}
0 & e^v   \\
-\bar{B}e^{-v} & v_{\bar{z}} - \frac{1}{2} H_{\bar{z}} e^{-v}
e^{\alpha}
\end{array}
\right) \left(
\begin{array}{c}
\psi_1 \\
\psi_2
\end{array}
\right)
\eeq
\end{lemma}

{\it Proof of Lemma.}
We have
$$
 \frac{\partial U}{\partial z} = v_z e^v = \frac{2H+i}{4}
  \psi_2 \partial \bar{\psi}_2
+ \frac {2H-i}{4} \bar{\psi}_1
  \partial \psi_1 - \frac{iH}{2} \psi_1 \bar{\psi}_2 |\psi_2|^2
  + \frac{H_z e^\alpha}{2}
$$
and combining that with (\ref{tildea}) we yield
$$
 \partial \psi_1 = (v_z - \frac{1}{2} H_z e^{-v} e^{\alpha}) \psi_1
  + \frac{1}{4} (2H+i) \widetilde{A} e^{-v} \psi_2,
$$
where $e^{\alpha} = |\psi_1|^2 + |\psi_2|^2 $.
Analogous calculations gives us
$$
\frac{\partial U}{\partial \bar{z}} = v_{\bar{z}} e^v =\frac{2H+i}{4}
\bar{\psi}_2 \bar{\partial} {\psi_2} +
\frac {2H-i}{4} \psi_1
\bar{\partial} \bar{\psi}_1 -
\frac{iH}{2} \psi_2 \bar{\psi}_1 |\psi_1|^2
  + \frac{H_{\bar{z}} e^\alpha}{2}
$$
and
$$
 \bar{\partial} \psi_2 = -\frac{1}{4} (2H-i) e^{-v} \bar{\widetilde{A}}
 \psi_1 + (v_{\bar {z}} - \frac{1}{2} H_{\bar{z}} e^{-v}
  e^{\alpha}) \psi_2.
$$
Together with the Dirac equation $\D\psi=0$ these equations constitute (\ref{dernil1}) and
(\ref{dernil2}). Lemma is proved.

Now let us prove the theorem. We again recall that  holomorphic
differentials on tori are constant: $\mathrm{const} \cdot dz^2$.
CMC surfaces in $\nil$ with $\widetilde{A}=0$ are spheres of
revolution \cite{AR2,BT2}. Hence $H$ and $\widetilde{A}$ are
nonvanishing constants and the equations (\ref{dernil1}) and
(\ref{dernil2}) are simplified as follows
$$
\bar{\partial} \left(
\begin{array}{c}
\psi_1 \\
\psi_2
\end{array}
\right) =
 \left(
\begin{array}{cc}
0 & e^v   \\
-\bar{B}e^{-v} & v_{\bar{z}} - \frac{1}{2} H_{\bar{z}} e^{-v}
e^{\alpha}
\end{array}
\right) \left(
\begin{array}{c}
\psi_1 \\
\psi_2
\end{array}
\right)
$$
which implies
$$
v_{z \bar{z}} + e^{2v} - |B|^2 e^{-2v}  = 0.
$$
By rescaling the conformal parameter we achieve that $|B|=1$. This proves
Theorem.

In this case the appearance of the same integrable system as the
Gauss--Codazzi equations for different classes of surfaces (CMC
tori in $\R^3$ and in $\nil$) does not mean any Lawson type
correspondence because for tori in $\nil$ this equation is written
not on the metric but on the potential $U$ of the Weierstrass
representation. \footnote{From the traditional point of view which
we do not follow, $U$ is not considered as a geometrical
quantity.} Moreover this coincidence does imply the local isometry
of corresponding surfaces.

We would like also to mention that until recently there are no
known examples of CMC tori in $\nil$ and this theorem is just a
step to proving their existence. One of the main difficulties is
that the systems (\ref{sinhgordon}) and (\ref{berdinsky}) are very
different from the physical point of view: they describe different
fields, i.e., the function $u$ in (\ref{sinhgordon}) is
real-valued and the function $v$ in (\ref{berdinsky}) in general
has nontrivial real and imaginary parts. Hence the reality
conditions for these systems are drastically different. However it
sounds possible to use soliton technique kind of the Lamb ansatz
to construct some analogs of the Abresch tori in $\R^3$
\cite{Abresch}.

{\bf 5. The spinor energy and the isoperimetric problem \cite{T2}}

Although in general for surfaces in $\nil$ and $\sll$ the
potentials $U$ and $V$ are complex-valued, the (spinor) energy
functional (\ref{energy}) is real-valued for compact oriented
surfaces without boundary. Moreover as in the Euclidean case it is
written in geometrical terms:

\begin{theorem}
[\cite{BT1}]
For a closed oriented surface $M$ in $G$ its (spinor) energy
$$
E(M) = \int_M UV dx dy
$$
equals
$$
\frac{1}{4} \int_M \left( H^2 + \frac{\widehat{K}}{4} -
\frac{1}{16} \right)
d \mu \ \ \ \ \mbox{for $G = \nil$};
$$
$$
\frac{1}{4} \int_M \left( H^2 +
\frac{5}{16} \widehat{K} - \frac{1}{4} \right) d \mu \ \ \ \
\mbox{for $G = \sll$,}
$$
where $\widehat{K}$ is the sectional curvature of the ambient space
along the tangent plane to the surface and $d\mu$ is the induced measure.
\end{theorem}

These expressions for $E$ are different from the Willmore functional
which for surfaces in a general ambient space is defined as
$$
{\cal W} = \int_M (|H|^2 + \widehat{K}) d\mu.
$$
For surfaces in $\R^3$ we have
$$
E = \frac{1}{4} {\cal W} = \frac{1}{4} \int_M
\left(\frac{\kappa_1+\kappa_2}{2}\right)^2 d\mu = \frac{1}{4}
\int_M \left(\frac{\kappa_1-\kappa_2}{2}\right)^2 d\mu +
\frac{1}{4} \int_M \kappa_1\kappa_2 d\mu.
$$
The Gauss--Bonnet theorem implies that for a compact oriented
surface $M$ without boundary we have
\beq
\label{willmore0}
E(M) = \frac{1}{4}\int_M \left(\frac{\kappa_1-\kappa_2}{2}\right)^2
d\mu + \frac{2\pi\chi(M)}{4}
\eeq
where $\chi(M)$ is the Euler characteristic of $M$. In particular this
implies that for spheres
$$
E \geq \pi
$$
and the equality is achieved exactly on the round spheres for which
$\kappa_1=\kappa_2$ everywhere.

We note that the round spheres are exactly the isoperimetric profiles in
$\R^3$, i.e. these are closed surfaces of minimal area among all surfaces
bounding domains of some fixed volume. It follows from
the variational principle that an isoperimetric profile is always
a CMC hypersurface at regular points and it is known that if the dimension of the ambient space
is not greater than seven then an isoperimetric profile is smooth.

The isoperimetric problem is not solved until recently for surfaces in
$\nil$. However it is known that in general for a compact Riemannian manifold
for small volumes the isoperimetric profiles are homeomorphic to a sphere
\cite{Morgan}. Hence for small volumes the isoperimetric profiles in $\nil$
are CMC spheres. By \cite{AR1} all CMC spheres are rotationally invariant, and by \cite{FMP},
CMC spheres of revolution form a family parameterized by
the mean curvature $H, 0 < H < \infty$.
We compute that

\begin{proposition}
For CMC spheres in $\nil$

\begin{enumerate}
\item
the energy functional is constant and equals
$E = \pi$;

\item the Willmore functional varies as follows:
$$
{\cal W}(H) = 10\pi  + \frac{\pi}{2H^2} -
$$
$$
- \pi\frac
{(1+4H^2)(3H^2-\frac{1}{4})H^3}{2} \left(\frac{\pi}{2}
 -\arctan\left [\frac{4H^2 -1}{4H}\right] \right).
$$
\end{enumerate}
\end{proposition}

Let us consider general surfaces of revolution in $\nil$. There is the natural
submersion
$$
\nil \to \nil/SO(2)
$$
onto the half-plane $u \geq 0$ with the metric
$$
du^2 + \frac{4dv^2}{4+u^2}.
$$
Let $\gamma(s) = (u(s),v(s))$ be a path-length parameterized
smooth curve in this halfplane which
generates by revolution a surface in $\nil$. Let us denote by $\sigma$
the angle between $\gamma$ and the vector
$\frac{\partial}{\partial u}$.
We have

\begin{theorem}
[\cite{BT2}]
For a closed oriented surface $M$ in $\nil$ obtained by revolving a curve
$\gamma \subset B$ around the $z$-axis, the spinor energy of $M$
equals
\begin{equation}
\begin{split}
\label{willmore1} E(M) = \frac 1 4  \int_{\gamma} \left( H^2 -
\frac 1 4 {n_3}^2 \right)
d\mu = \\
\frac{\pi}{8} \int_{\gamma} \left(\dot {\sigma} -
\frac{\sin \sigma}{ u}\right)^2 \sqrt{4u^2 + u^4} ds - \frac{\pi}{4}
\int_{\gamma} \frac{\partial{[\dot{u}
\sqrt{4+u^2}]}}{\partial s} ds = \\
\\
\frac{\pi}{8} \int_{\gamma} \left(\dot {\sigma} -
\frac{\sin \sigma}{ u}\right)^2 \sqrt{4u^2 + u^4} ds + \frac{\pi \chi(M)}{2}
\end{split}
\end{equation}
where $\chi(M)$ is the Euler characteristic of $M$.

If $\dot{\sigma} = \frac{\sin \sigma}{u}$ everywhere then
the surface is a CMC sphere.
\end{theorem}

It implies

\begin{corollary}
For spheres of revolution in $\nil$ we have
$$
E(M) \geq \pi
$$
and the equality is
attained exactly at CMC spheres.
\end{corollary}

\begin{corollary}
For tori of revolution in $\nil$ the spinor energy is positive:
$$
E(M) > 0.
$$
\end{corollary}

It is also straightforward to prove

\begin{proposition}
[\cite{BT2}]
The CMC spheres in $\nil$ are the critical points of the spinor
energy  functional $E$.
\end{proposition}

We see now that except the spectral theory of the Dirac operator
there are other reasons to treat the spinor energy as the right analog of the
Willmore functional for surfaces in $\nil$. Indeed,

\begin{itemize}
\item
it takes the constant value on the CMC spheres which are the critical points of
this functional;

\item
there is a strong similarity of formulas (\ref{willmore0}) and
(\ref{willmore1}). However the quantities
$\dot{\sigma}$ and $\frac{\sin \sigma}{u}$
are not the principal curvatures of a surface of revolution
and two poles are the only umbilic points on a CMC sphere in $\nil$;

\item
the conditions $A=0$ and $\widetilde{A}=0$
distinguish in $\R^3$ and $\nil$ the minima of $E$ for spheres of revolution
(in the Euclidean case even for spheres).
\end{itemize}

Of course, this study has to be completed and the following questions
are worth to be answered:

\begin{enumerate}

\item
{\sl is $E$ bounded from below for each topological type
of closed oriented surfaces?}

\item
{\sl is  $E$ positive?}

\item
{\sl are the CMC spheres in $\nil$ are the global minima of $E$ for spheres?}

\item
{\sl how to generalize \eqref{willmore1} for general surfaces?}

\item
{\sl what are the minima of $E$ for surfaces of fixed topological type
and, in particular, what is the substitution of the Willmore
conjecture?}
\end{enumerate}

It is also interesting to study the analogous questions for surfaces in $\sll$
for which the spinor energy functional also has a geometrical form.

For $S^2 \times \R$ we have the following computational observation:

\begin{proposition}
[\cite{BT2}]
For isoperimetric profiles $M$ in $S^2 \times \R$ we have
$$
\int_M (H^2 + \widehat{K} +1)d\mu = 16\pi.
$$
\end{proposition}

The isoperimetric problem for $S^2 \times \R$ was solved by Pedrosa
\cite{Pedrosa} who proved that for volumes $d \leq d_0$ the isoperimetric
profiles are CMC spheres, for $d>d_0$ the isoperimetric profiles bound
the product cylinders $S^2 \times \left[0,\frac{d}{4\pi}\right]$ where
$d_0$ is some transition point from one topological class of solutions
to another. The functional mentioned in Proposition takes the same value
on all CMC spheres (not only isoperimetric) and on all isopermetric profiles
(connected and disconnected).

We would like to guess that

{\it the right analog of the Willmore theory (at least for spheres) has to be related to the
isoperimetric problem and the isoperimetric profiles in three-dimensional
homogeneous spaces have to be distinguished as (at least local)
minima of the Willmore type functional which is constant on them.}

\end{document}